\documentstyle{amsppt}
\magnification=1000
\def\ss{\smallskip}

\def\oonL{\mathop{{\hbox{$L^{m,2}$}\kern -20pt\raise7pt
\hbox{$\circ$}}}}

\def\O{\Omega}

\def\bigtimes{\mathop{{\lower2.95pt\hbox{$\wedge$}\kern-6.666pt\raise2.95pt
\hbox{$\vee$}}}\limits}


\def\R{${\text{\bf R}}^N$}
\def\no{\noindent}

\def\n{\nabla}

\def\da{d_{\partial \O}(x)}
\def\d{\da}
\def\pn{\par\noindent}

\let\hacek=\v
\def\v{\vert}

\hsize=5in
\vsize=7.6in
\baselineskip=16pt

\def\Wmps0O{W^{m,p}_0(\O,\d^s)}
\def\Wmpt0O{W^{m,p}_0(\O,\d^t)}

\nologo
\topmatter
\title
A remark on the history of Hardy inequalities in domains.
\endtitle
\author
Andreas Wannebo
\endauthor
\endtopmatter

\document

\no
This is a text on the history of Hardy inequalities in domains that are
subsets of \R.
Hardy inequalities for domains are a generalisation of the one-dimensional
Hardy inequality. In one dimension there is no geometrical problem
coming from the domain. 
For Hardy inequalities in domains the situation is much more difficult. 
\ss
\no
To begin with we list some of the limitations of this account.
The history of Hardy inequalities (in domains) includes several aspects.
The following is left out:
Norms other than $L^p$-norms (weights),
trivial domains in $R^N$,
general manifolds,
the best constants aspects and the cases of
results that so far have been unproductive with respect to Hardy inequalities.
-- The goal is to present the main line of progress in the area over the
years. However with no claim of completeness.
\ss
\pn
The following formula describes a rather general Hardy 
inequality in a domain
$$
(\int_{\O}|u|^q\d^bdx)^{1\over q} \le A (\int_{\O}|\n^mu|^p\d^adx)^{1\over p}.
$$
The notation is as follows.
The domain in Euclidean N-space is denoted $\O$,
the distance function from a point inside the domain to the boundary 
is denoted by $\d$,
the function $u$ is in say $C^{\infty}_0(\O)$, i.e. infinitely differentiably
functions with compact support in $\O$,
$m$ is a positive integer and the $m$th gradient is denoted $\n^mu$.
This gives power type formulae, which are the most important ones.
\ss
\pn
The study of Hardy inequalities for domains got its beginning with
studies by J.Ne\hacek cas [1]. He carried over the 
one-dimensional power-type Hardy inequality to bounded Lipschitz domains. 
Since this is a kind of test case regarding the parameter space
we give a more complete statement.
The theorem states that,  
if given $m=1$ (i.e. a first order case), 
$1<p=q<\infty$, 
$a<p-1$ and $b=a-p$,
then for each bounded $\O$ of Lipschitz domain, there exists a constant $A$ 
such that for all functions in $C^{\infty}_0(\O)$ the inequality formula
above holds. This is a dilation homogeneous formulation.
Furthermore the proof gives the value of a possible 
constant.
-- Observe that the possible exponents here are exactly the same as in the 
one-dimensional case.
The method of proof can be named ``the one-dimensional method''.
It is based on a study of the function on sets of lines covering the domain.
The final result arises from integration over inequalities on these lines.
\ss
\no
In the study of Hardy inequalities for domains it seems, unlike the 
one-dimensional case, not possible to get complete results 
(i.e. {\bf results}) for all weights.
Hence the study is more concentrated to smaller weight classes.
The most important ones are power type weights or weights that are 
functions of the distance to boundary. 
\pn
-- From a functional analysis or potential theoretic point of view 
PDE questions can be treated with help of Hardy inequalities.
Hence an improvement of the results for Hardy inequalities by an enlargement
of the class of domains satisfying a certain fixed Hardy inequality 
implies that often a larger class of domains can be treated with similar 
arguments and results.
\ss
\no
In the middle-late 60ies A.Kufner develops the one-dimensional method
in order to treat ``bad'' domains. Here
``bad'' is H\"older domains. To understand what these domains look like
a rough idea is to take a domain which is defined to have a power
type cusp with fixed power that can be placed inside the domain
with the endpoint touching any given boundary point. 
Much of his results are presented in ``Weighted Sobolev Spaces'', [2]. 
The idea is that corresponding to the ``badness'' of the domain
new sets of parameters are used in the Hardy inequality formula.
This implies a corresponding weakening of the inequality.
The condition on the parameters are then dependant on
the H\"older class parameter, cf. the power of the cusp described.
Kufner also treats a situation where there is no boundary condition.
The inequality is then saved by adding a natural term on the 
right hand side. The result is a natural companion to the previous case.
A situation treated by Kufner is that with boundary condition only on
some subset of the boundary. 
There is also a treatment of more general weights in [2].
\ss
\no
There are many papers on Hardy inequalities in domains by many authors
following the first papers by Kufner. However as we shall see
the development later took a different turn and this calls for 
a jump in the story told here. 
However let us just mention a work done by J.Kadlec and Kufner [3] and [4]
as of special interest.
\ss
\no
A.Ancona in 1981 [5] published a short note in order to give substance to
the results planned to appear in a paper by F.Browder and H.Brezis [6]. 
The note included two
results. The last was a Hardy inequality of order one where the requirement
was that $p>N$ (and boundedness which happens to be unnnecessary).
Otherwise the parameters are as in the Ne\hacek cas case above but with
weights restricted by $a=0$. This is a surprising result destroying old ideas
of ``bad'' and ``good'' domains.
The first result concerned the question when the Sobolev space determined
by $C^{\infty}_0(\Omega)$ has its functions given as difference
of two nonnegative such functions.
\ss
\no
The present author, then student of L.I.Hedberg, was given the assignment to 
study the Ancona note. This in fact unexpectedly lead to a 
prolonged study of Hardy inequalities but now with entirely new ideas.
These involved Poincar\'e inequalities with a closed subset 
determining the zeros of the functions on a ``nice'' domain and covering 
of the ``bad'' domain with scaled copies of the ``nice'' domain with
centres at the boundary of the ``bad'' one. This causes an infinite
overlap which was handled by a new summation method.
At this point the work of V.G.Maz'ya came into the picture
since these Poincar\'e inequalities for a ``nice'' domains 
could be handled by a result by him
on Poincar\'e inequalities for a cube 
formulated in terms of polynomial capacities,
invented by him for this purpose. 
His best account of this is one chapter of his book ``Sobolev Spaces'' [7].
These polynomial capacities can readily be calculated or estimated
i.e. they are ``nice'' despite a somewhat ugly definition. The order
of the polynomial capacity 
is the order of the gradient in the corresponding Poincar\'e inequality.
If the order is one then the polynomial capacities are just
ordinary (non-linear) capacities from (non-linear) potential theory.
When the order is higher the polynomial capacities becomes bigger
and then the Poincar\'e inequality gets the same or better constant 
-- this means smaller. This can give new domains for the related 
Hardy inequalities.
The polynomial character of polynomial capacities appears first with
higher order. Let the order be $m$, then the polynomials in the 
polynomial capacities are the null set of the $m$th gradient 
(as an operator),  i.e. they are the polynomials of degree $\le m-1$. Take the 
order two case, then the null set consists
of the polynomials of degree one or zero.
One aspect now is how a given point set in a cube
say, deviates from every one of the zero manifolds of degree one 
polynomials, i.e. all hyperplanes. 
The degree of approximation here is measured by
order two non-linear capacities (nicer i.e. bigger than order one non-linear
capacities). 
In the degree zero case the 
approximation is measured in some sense by order one non-linear capacities.
-- To sum up, the construction led to a sufficient condition for
Hardy inequalities with the weight powers as in the Lipschitz case above.
The condition for suffiency was what later by Lewis called
uniform p-fatness of the complement set.
In the order one case this can be stated as uniform fatness in terms
of order one capacities and but in the higher order cases 
the situation is better and is formulated in terms
of uniform fatness of a the higher order polynomial capacity.
At this time of the story the present author only had proved this 
for $a<0$, $q=p$ and $m$ general. The interesting fact was that
the inequality was proved for different $a<0$ with the same conditions.
The problem was that the constant $A$ went to infinity when $a$ went to 
zero. The key question at the time.
-- At this point the present author and Ancona met in Uppsala 
and the situation was explained to him.
In order to get the simplest case of annoyance consider general domains 
in $\text{\bf R}^2$ with complement satisfying the uniform 2-fat condition
and with $m=1$, $p=q=2$ and $a=0$.
\ss
\no
Ancona solved the problem, see [8]. 
This involved $m=1$, $p=q=2$, $a=0$ and $b=-2$ for dimension $N$.
A striking example here is the following, 
given any simply connected domain in the plane, not the plane itself, 
then the above Hardy inequality holds for the parameters given.
Here a standard projection property of capacities has been used.
He proved that generally the uniform 2-fatness for the complement
is sufficient with these parameters, but for $N=2$ this is  necessary as well.
There are other results as well in [8].
The paper has had influence, being cited 28 times so far. 
Ancona used that he got the question together with 
the correct formulation from the present author and by a mistake 
this was not mentioned. We have had plenty of time afterwards though
to talk this over and this is not any matter of disagreement.
\ss
\no
J.L.Lewis was referee of [8], that way he got interested in the question.
He was able in, see [9], to enlarge on the Ancona result.
His main results are that the Hardy inequality holds for $m=1$, $p=q>1$,
$b=a-p$ and $a\le 0$, if the domain has a uniformly $p$-fat complement.
Morever if $p=N$,
then $p$-fatness of the complement also is necessary.
\ss
\no
A remark. The first result in the 
Ancona note [5] -- read carefully, 
together with the Ancona paper [8] or Lewis [9] has the following implication.
Let $W^{1,p}(\O)$ the Sobolev space of order one with domain $\O$ and 
parameter $p$. If the domain has a complement which is uniformly 
$p$-fat,
then $W^{1,p}(\O)=W^{1,p}_0(\O)$, i.e. it equals the Sobolev space 
induced by $C^{\infty}_0(\O)$, if the Hardy inequality holds for
$m=1$, $p=q$, $a=0$ and $b=-p$.
\ss
\no
From the middle 80ies and on the present author made several
typed respectively TEXed notes. (Finally more then 50pp.)
They were circulated. The first publication was [10], 1990.
It was written because at this time the present author,
who was not aware of the work of Ancona [8] and Lewis [9],
had overcome this $a=0$ problem. The paper was written
in order to document a small portion of these the notes.
(Reference is made.)
The fundamental summation method here is not the one in these
notes however, 
but a rather equivalent one.
The contents of [10] is 
generalizing the results
of Ancona and Lewis. The range of weights is somewhat larger,
$a<s_0$, with $0<s_0$ and higher order inequalities are 
treated on equal basis. Also it is shown how these 
polynomial capacities get into the picture and reference is given 
to some of their properties.
The case $p=1$ is not explicitly treated, but by only reading the beginning
of the proof it is easily seen that in this case a $log$ 
factor enters in the right hand side where otherwise the parameters
are $m$, $q=q=1$, $b=-m$, and $a=0$.
\ss
\no
At last the present author made two preprints  [11] and [12] in 
order to cover most of the material in these circulated notes. 
In [11] the development of polynomial capacities is taken further
both conceptually and given more results. One application is made 
to the description of
certain subspaces and subcones of Sobolev space.
In [12] this new knowledge of polynomial capacities is used to make
a better treatment of Hardy inequalities. Furthermore 
extra ideas concerning the construction of the Hardy inequalities is given.
These results are then given both generally as well as expressed
in geometric etc. terms.
There are far too many results here and hence a solution has been 
to present them in a catalogue kind of way.
The original question by Ancona [5], i.e. when holds the following for
Sobolev space and its nonegative cones
$$
W^{m,p}_0(\O)=W^{m,p}_0(\O)_+-W^{m,p}_0(\O)_+,
$$
is given answers of quite general kind, etc. A strong conjecture is made.
\ss
\no
{\bf Conjecture:} The decomposition above holds for all $m$ odd,
$p>1$, $N$ and open $\O\subset {\text \R}$ and for 
$m$ even there are $p$, $N$, $\O$ such that it is untrue.
\ss
\pn
Further preprints are planned to follow up [10] and [11].
\ss
\pn 
Remark: By calculations the present author has shown 
that if the question of Hardy inequalities is specialized to order one 
and H\"older domains then there is a part of 
the region of parameter values where the old parameter
values given by Kufner are better as description of the Hardy 
inequalities than what is got by the polynomial capacity
approach in [11].
\ss
\no 
The present author made a theorem in the thesis 1991, see [13], which
was intended to test how far the one-dimensional method could go.
It contains all the Hardy inequalities in the Kufner 
book (except those which only have zero on part of the boundary,
which were not tried).
This gives much wider range of domains and weights. 
\ss
\no
Recently there have been established a connection between pointwise 
maximal inequalities and Hardy inequalities which have caught interest. 
Here are some of these papers are listed in order after dates.
They are given here without further evaluation,
P.Hajlaz: ``Pointwise Hardy inequalities'' [14];
J.Kinnunen and O.Martio: ``Hardy's Inequalities for Sobolev Functions'' [15];
D.E. Edmunds and J. R\'akosn\'ik: ``On a Higher-Order Hardy Inequality'' [16].


\Refs

\ref
\no 1
\by J.Ne\hacek cas
\paper Sur une m\'ethod pou r\'esondre les \'equations aux d\'eriv\'ees
partielle du type elliptique, voisine de la variationell.
\jour Ann. Scoula Norm. Sup. Pisa
\vol 16
\yr 1962
\pages 305-326
\endref

\ref
\no 2
\by A.Kufner
\book
Weighted Sobolev Spaces
\yr 1980
\publ Teubner-Texte zur Mathematik, Leipzig
\endref

\ref
\no 3
\by Kadlec and A.Kufner
\paper Charcterization of functions with zero traces by integrals
with weght functions I.
\jour \hacek Casopis P\hacek e. Mat.
\vol 91
\yr 1996
\pages 463-471
\endref
 
\ref
\no 4
\by Kadlec and A.Kufner
\paper Charcterization of functions with zero traces by integrals
with weght functions II.
\jour \hacek Casopis P\hacek e. Mat.
\vol 92
\yr 1997
\pages 16-28
\endref

\ref
\no 5
\by A.Ancona 
\paper Une propri\'et\'e des espaces de Sobolev
\pages 477-480
\jour C.R.Acad. Sc. Paris
\vol 292 
\yr 1981
\endref

\ref
\no 6
\by H.Brezis and F.Browder
\paper Some properties of higher order Sobolev spaces
\pages 245-259
\jour J. Math. Pure. Appl.
\vol 61
\yr 1982
\endref

\ref
\no 7
\by V.G.Maz'ja
\book Sobolev spaces
\publ Springer Verlag
\yr 1985
\endref

\ref
\no 8
\by A.Ancona
\paper On strong barriers and the inequality of Hardy for domains in $R^n$
\pages 274-290
\jour J. London Math. Soc. 
\vol 34 
\yr 1986
\endref

\ref
\no 9
\by J.L.Lewis
\paper Uniformly fat sets 
\pages 177-196 
\jour Trans. Amer. Math. Soc. 
\vol 308 
\yr 1988
\endref

\ref
\no 10
\by A.Wannebo
\paper Hardy inequalities 
\pages 85-95
\jour Proc. Amer. Soc. 
\vol 109 
\yr 1990
\endref

\ref
\no 11
\by A.Wannebo
\paper Polynomial capacities, Poincar\'e type inequaities and 
spectral analysis in Sobolev space
\pages 1-33
\jour Trita-Mat-1998-42
\vol 
\yr 1998
\paperinfo{Report series, Math, KTH, Sweden}
\endref

\ref
\no 12
\by A.Wannebo
\paper Hardy and Hardy PDO type inequalities in domains part I
\pages 1-31
\jour Trita-Mat-1999-15
\vol 
\yr 1999
\paperinfo{Report series, Math, KTH, Sweden}
\endref

\ref
\no 13
\by A.Wannebo
\paper Hardy inequalities and imbeddings in domains generalizing $C^{0,\lambda}$ domains. 
\pages 1181-1190
\jour Proc. Amer. Soc. 
\vol 122 
\yr 1994
\paperinfo{Part of thesis 1991}
\endref

\ref
\no 14
\by P.Hajlasz 
\paper Pointwise Hardy inequalities
\pages 417-423.
\jour Proc. Amer. Soc. 
\vol 127, 2
\yr 1999
\endref

\ref
\no 15
\by J.Kinnunen and O.Martio
\paper  Hardy's inequalities for Sobolev functions
\pages 489-500
\jour Math. Res. Lett. 
\vol 4 
\yr 1997
\endref

\ref
\no 16
\by D.E.Edmunds and J. R\'akosn\'ik
\paper On a higher order Hardy inequality
\pages 113-121
\jour Mat. Bohem.
\vol 124, 2 
\yr 1999
\endref

\endRefs
\enddocument